\documentclass[a4paper, 11pt]{article}                     

\usepackage{graphicx}
\usepackage{url}
\usepackage{hyperref}

\usepackage{amssymb}

\usepackage{color, colortbl}

\usepackage{url}  
\usepackage{ifthen}  
\usepackage{multicol}   
\usepackage{multirow}
\usepackage[utf8]{inputenc} 
\usepackage{rotating}
\usepackage{algorithm}
\usepackage{algorithmic}
\usepackage{comment}
\usepackage{amsmath}
\usepackage{mathtools}
\usepackage{graphics}
\usepackage{array}
\newcolumntype{C}[1]{>{\centering\let\newline\\\arraybackslash\hspace{0pt}}m{#1}}

\usepackage{threeparttable}
\usepackage{mathrsfs}
\usepackage{cprotect}
\usepackage{arydshln}
\usepackage{ dsfont }
 \usepackage{amsthm}
 \usepackage{comment}
\usepackage[normalem]{ulem}
\usepackage{bm}
\usepackage{amsfonts}

\newtheorem{theorem}{Theorem}
\newtheorem{remark}{Remark}
\newtheorem{proposition}{Proposition}

\newcommand{\var}{\mbox{var}}

\usepackage[round]{natbib}

 \usepackage{mathptmx}      

\usepackage{latexsym}

\newcommand{\institute}[1]{#1}
\newcommand{\email}[1]{\href{mailto:#1}{#1}}
\newcommand{\keywords}[1]{\textbf{ \textsc{Keywords: #1}}}

\begin{document}

\title{A Permutation Approach to Goodness-of-fit Testing in Linear Regression Models
}


\author{Rok Blagus         \and
        Jakob Peterlin \and
        Janez Stare 
}


        \date{ 2019-11-15}

\maketitle

\textit{\footnotesize{
\institute{Rok Blagus
              Institute for Biostatistics and Medical Informatics, University of Ljubljana\\
              Vrazov trg 2, Ljubljana, Slovenia \\
              \email{rok.blagus@mf.uni-lj.si}
              }
}}

\begin{abstract}
 Model checking plays an important role in linear regression as model misspecification seriously affects the validity and efficiency of regression analysis. In practice, model checking is often performed by subjectively evaluating the plot of the model's residuals. This approach is objectified by constructing a random process from the model's residuals, however due to a very complex covariance function obtaining the exact distribution of the test statistic is intractable. Several solutions to overcome this have been proposed, however the simulation and bootstrap based approaches are only asymptotically valid and can, with a limited sample size, yield tests which have inappropriate size. We therefore propose to estimate the null distribution by using permutations. We show, under some mild assumptions, that with homoscedastic random errors this yields consistent tests under the null and the alternative hypotheses. Small sample properties of the proposed tests are studied in an extensive Monte Carlo simulation study, where it is demonstrated that the proposed tests attain correct size, even with strongly non-normal random errors and a very small sample size, while being as powerful as the other available alternatives. The results are also illustrated on some real data examples.

\end{abstract}

\section{Introduction}

Linear regression is a fundamental statistical tool for analyzing 
data for a continuous outcome. Model misspecification seriously affects the validity and efficiency of regression analysis, therefore model checking plays an important role. In this paper we propose goodness-of-fit tests for the linear regression model, where, given a random sample $(Y_i,x_{i1},...,x_{i(p-1)})$, $i=1,...,n$, the relation between the outcome, $Y_i$ and the assumed fixed, possibly transformed, covariates, $x_{ij}$, $j=1,...,p-1$ is given by

\begin{equation}
\label{eq1.1}
Y_i=\beta_0+\beta_1x_{i1}+...+\beta_px_{i(p-1)}+\epsilon_i,
\end{equation}
where $\epsilon_i$ are independent and identically distributed (i.i.d.) random variables with expectation zero, $E(\epsilon_i)=0$, and a constant variance, $\var(\epsilon_i)=\sigma^2<\infty$, $\beta_0$ is the intercept and $\beta_j$, $j=1,...,p-1$ are regression coefficients. The tests are based on residuals, $e_i=Y_i-\hat Y_i$, where $\hat Y_i=\hat\beta_0+\hat\beta_1x_{i1}+...+\hat\beta_{p-1}x_{i(p-1)}$ are fitted values for subject $i$ based on 
estimates $\hat\beta_j$. The residuals are standardized and ordered by their respective fitted values and a stochastic process,  based on a cumulative sum of (standardized and ordered) residuals, \begin{equation}
\label{eq1}
W_n(t)=\frac{1}{\sqrt{n\hat\sigma^2}}\sum_{i=1}^ne_i I(\hat Y_i \leq t), t\in \mathds{R},
\end{equation}
where $\hat\sigma^2$ is a consistent estimator of $\sigma^2$ and $I(\cdot)$ is the indicator function, is constructed.

Test statistic, $T$, used for testing the null hypothesis,
$$H_0: \mbox{model (\ref{eq1.1}) generated the data} $$
can then be any function of $W_n(t)$ such that large (or small) values correspond to evidence against $H_0$. 
The tests formalize the graphical procedures where the residuals are plotted versus the fitted values where one looks for systematic patterns. Under the null hypothesis there is no systematic pattern and the residuals are centered around zero. The appearance of a systematic pattern may indicate a lack-of-fit. 
The proposed tests will have good power against the alternatives which cause a systematic pattern in terms of the successive number of positive versus negative ordered residuals.

Due to dependence between the residuals exact distribution of even the most common test statistics (e.g. Kolmogorov-Smirnov (KS) type test statistic and Cramer-von Mieses (CvM) type test statistic) is unknown. Exact results which are available for the KS and CvM test statistic for the Brownian motion and Brownian bridge stochastic processes rely on the assumption of independence \citep{Beghin99,Kulperger90}, hence they would, due to negative correlation of the residuals, yield conservative tests. 
For example, \citet{Royston92} proposed a goodness-of-fit test, which is based on the similar process as (\ref{eq1}). However, using the standard normal distribution to obtain $p$-values, thus neglecting the dependence of the residuals, led to extremely conservative test where the null hypothesis was never rejected (see \citet{HosmerStatMed}).

We propose to overcome this issue by using permutations to approximate the null distribution of the test statistic. Let $T^k$ be a test statistic in permutation $k$, $k=1,...,K$, where $K$ is either the number of all possible permutations or a large number of random permutations \citep{dwass57,phipsonsmith10,Hemerik2018}. The $p$-value is then defined in the usual way as the proportion of $T^k$ which are as or more extreme as $T$.



Using permutations with linear models is common but has not been applied in the context of goodness-of-fit testing (see \citet{diCiccio16} and references therein). There are two different permutation procedures which one could adapt to this setting. One is raw data permutation procedure proposed by  \citet{Manly06} and the other is permutation of residuals proposed by \citet{Braak92}, which in our setting is equivalent to the permutation procedure proposed by \citet{Friedman83}. Both procedures were shown to provide valid inference in the context of tests for partial regression coefficients in a linear model \citep{Anderson99,Anderson01,diCiccio16}. In this paper we show, under some mild assumptions, that, conditional on the data, the permutation of residuals, under $H_0$ well approximates the distribution of the constructed random process and therefore provides valid inference in the context of goodness-of-fit testing.

Small sample properties of the proposed tests are studied by an extensive Monte-Carlo simulation study.  
It is shown that our proposed tests perform well also when the errors are non-normal and the sample size is limited. The tests are illustrated also on data from a study of fetal mandible length \citep{Chitty93} and data from \citet{draper81}.




\subsection{Existing approaches}

Several goodness-of-fit tests based on a similar, cumulative sums process, as defined in (\ref{eq1}) 
were proposed \citep{Lin02,Stute98a,Diebolt01,Su91,Fan01,Stute98b}. 
They are mainly defined as
\begin{equation}
\label{eq0}
\tilde W_n(\mathbf{t})=\frac{1}{\sqrt{n}}\sum_{i=1}^ne_i I(\mathbf{x}_i\leq \mathbf{t}),
\end{equation}
$\mathbf{x}_i=(1,x_{i1},...,x_{ip})^T\in \mathds{R}^{p}$, $\mathbf{t}=(t_0,t_1,...,t_{p})^T\in \mathds{R}^{p}$, $I(\mathbf{x}_i\leq \mathbf{t})=I(x_{i1}\leq t_1,...,x_{ip}\leq t_{p})^T$, where $\mathbf{c}^T$ denotes the transpose of $\mathbf{c}$, and the test statistics are then KS or CvM statistics based on $\tilde W_n(\mathbf{t})$.  

Note that in contrast to our proposed test, the residuals in (\ref{eq0}) are not standardized which can decrease the power of the tests \citep{Christensen15} and can under $H_0$ slow the convergence rate. Also, in (\ref{eq0}) multivariate ordering procedure is used which can be problematic in higher-dimensions \citep{Christensen15}, 
therefore \citet{Lin02} also proposed to order the residuals in (\ref{eq0}) by the fitted values. They also considered to take the sum only within a window specified by some positive constant $b$, i.e. in (\ref{eq0}) they use $I(t-b<x_{ij}\leq t)$. 
This is done as the process $\tilde W_n(t)$ tends to be dominated by the residuals with small covariate values. All these different ways of ordering the residuals could easily be used also in our proposed test, but we believe that ordering the residuals by the fitted values is the most natural option when one wants to test for the possible lack-of-fit of the entire fitted model. We show later how some more novel ways of ordering the residuals can be defined to obtain tests which target the lack-of-fit of different parts of the model.  

The existing tests differ in the way how the null distribution of the proposed test statistic is obtained. The result of \citet{Stute97} depended on the asymptotic distribution of $\tilde W_n(\mathbf{t})$ but did not yield satisfactory results with small sample size \citep{Stute98a}. \citet{Su91} and \citet{Lin02} used a simulation approach to obtain the $p$-values which relies on approximating $\tilde W_n(\mathbf{t})$ with a random process which for the linear model is defined as,
\begin{equation}
\label{sw}
\hat W_n(\mathbf{t})=\frac{1}{\sqrt{n}}\sum_{i=1}^n \left[ I(\mathbf{x}_i\leq \mathbf{t})-\left(\sum_{k=1}^n   \mathbf{x}_k^T I(\mathbf{x}_k\leq \mathbf{t}) \right)\left[\sum_{k=1}^n \mathbf{x}_k\mathbf{x}_k^T \right]^{-1}\mathbf{x}_i  \right]e_iZ_i,
\end{equation}
where $(Z_1,...,Z_n)$ are independent standard normal variables (see \citet{Su91} and \citet{Lin02} for more details). Note that $\hat W_n(\mathbf{t})$ can equivalently be defined as
 $$\hat W_n(\mathbf{t})=\frac{1}{\sqrt{n}}\sum_{i=1}^n e_i^*I(\mathbf{x}_i\leq \mathbf{t}),$$
where $e_1^*,\ldots,e_n^*$ are the residuals obtained when regressing $\mathbf{Y}^*$ on $\mathbf{x}$, where $\mathbf{Y}^*=(Y_1^*,\ldots,Y_n^*)^T$ and
\begin{equation}\label{ee1}Y_i^*=\hat{Y}_i+e_iZ_i\mbox{, }i=1,\ldots,n.\end{equation}
 This approach is however only asymptotically valid \citep{Stute98a} and, as it is shown later, with small sample size yields tests which do not attain the nominal level. \citet{Stute98a} used bootstrap \citep{Wu86,Hardle93,Liu88} to obtain $e_1^*,\ldots,e_n^*$. They considered classical, smooth, wild and residual bootstrap. In wild bootstrap, $Z_i$ in (\ref{ee1}) is replaced by $V_i$, where $(V_1,\ldots,V_n)$ are independent and identically distributed such that $E(V_i)=0$, $var(V_i)=1$ and $|V_i|\leq c<\infty$ for some finite $c$. In residual bootstrap, $e_iZ_i$ in (\ref{ee1}) are replaced by an iid sample from the empirical distribution function of the (centered) $e_i$'s. 

Recently, \citet{Hattab18} proposed several tests based on partial sums of residuals where the test statistics are based on sums of a subset of the (ordered and standardized) residuals \citep{Christensen10,Christensen15}. Defining partial sums over a subset of the residuals, as opposed to cumulative sums over all residuals as in (\ref{eq1}), makes it possible to determine the asymptotic distribution,  when the (appropriately standardized) test statistic is based on the largest of the absolute values of partial sums \citep{Christensen15} or the largest partial sum of absolute values of the (ordered and standardized) residuals \citep{Hattab18}. The rate of convergence is however very slow and using only a subset of the residuals can lead to loss of power when ordering is not appropriately chosen \citep{Hattab18}. Determining the number of residuals which are included in the subset for which the test statistic is computed can also be problematic in practice. To solve the issue with a slow rate of convergence, \citet{Hattab18} proposed a Monte-Carlo simulation scheme. In the case when the partial sums are based on all residuals, this approach is identical as the approach proposed by \citet{Su91} and \citet{Lin02}. While the tests based on partial sums can, if the residuals are ordered appropriately, under some alternatives be more powerful than the tests based on cumulative sums, they are not considered here in more detail, since they lack a nice visual representation which is one of the main strengths of the approach considered here \citep{Lin02}. The approach considered here could however easily be used instead of the Monte-Carlo procedure used by \citet{Hattab18} also for tests based on partial model checks.

In the Supplementary material 2 we show, using the same simulation setup as is used in \citet{Stute98a} and \citet{Christensen15}, that our proposed tests attain correct nominal level and are at least as powerful.

\section{The proposed goodness-of-fit tests and their asymptotic convergence}
Here we first introduce some additional notation and formally present the proposed goodness of fit test. This is then followed by an asymptotic study of the convergence of the constructed random processes under the null hypothesis. Finally, the asymptotic convergence of the constructed random processes under a particular alternative hypothesis is studied in the last section. The proofs of the theoretical results are presented in Supplementary material 1.

\subsection{Notation}

Write the model in the matrix form as
\begin{equation*}\bm{Y}=\mathbf{x}\bm{\theta}+\bm{\epsilon},\end{equation*}
where $\bm{Y}=(Y_1,\ldots,Y_n)^T$ and $\bm{\epsilon}=(\epsilon_1\ldots,\epsilon_n)^T$ are $n$-vectors of response variables and errors, respectively, $\bm{\theta}=(\beta_0,\beta_1,...,\beta_{p-1})^T$ is a $p$-vector of model parameters and $\mathbf{x}\in \mathds{R}^{n\times p}$ is assumed fixed real design matrix of full rank where the first column of $\mathbf{x}$ is a $n$-vector of ones, $\bm{1}=(1,\ldots,1)^T$. Let $\mathbf{x}_i$, a $p$-vector, be the $i$-th row of the design matrix. With $\mathbf{P}=\mathbf{x}(\mathbf{x}^T\mathbf{x})^{-1}\mathbf{x}^T$ we denote the projection matrix. The residuals are then
$$\bm{e}=(e_1,...,e_n)^T=\bm{Y}-\bm{\hat Y}=\bm{Y}-\mathbf{x}\bm{\hat\theta}=(\mathbf{I}-\mathbf{P})\bm{Y},$$
 where
 $$\bm{\hat\theta}=(\mathbf{x}^T\mathbf{x})^{-1}\mathbf{x}^T\bm{Y}$$
 is an ordinary least-squares (OLS) estimator of $\bm{\theta}$.

 Denote the OLS estimator and residuals obtained when regressing $\bm{Y}^k$ on $\mathbf{x}$ as
 $$\bm{\hat\theta}^k=(\mathbf{x}^T\mathbf{x})^{-1}\mathbf{x}^T\bm{Y}^k$$
 and
 $$\bm{e}^k=(e_1^k,...,e_n^k)^T=\bm{Y}^k-\bm{\hat Y}^k=\bm{Y}^k-\mathbf{x}\bm{\hat\theta}^k=(\mathbf{I}-\mathbf{P})\bm{Y}^k=(\mathbf{I}-\mathbf{P})\bm{\Pi}\bm{e},$$
 respectively, where
$$\bm{Y}^k=(Y_1^k,...,Y_n^k)^T= \mathbf{x}\bm{\hat\theta}+\bm{\Pi}\bm{e},  $$
and $\bm{\Pi}\in \mathds{R}^{n\times n}$ is some permutation matrix. Obviously, when $\bm{\Pi}=\mathbf{I}$, $\bm{e}^k=\bm{e}$. 

\subsection{The proposed goodness-of-fit tests}

After the model has been estimated and the model's residuals, $e_i$, have been obtained, these are then ordered by their respective fitted values, $\hat Y_i$,  to obtain the vector of ordered residuals, which is then standardized by dividing with 

$$ \hat\sigma=\sqrt{\frac{1}{n-p}\sum_{i=1}^n e_i^2}. $$


(Any other consistent estimator of $\sigma$ could be used instead of $\hat\sigma$.) A stochastic process is then constructed from the vector of ordered standardized residuals by taking the cumulative sum of the vector's elements, obtaining the process defined in (\ref{eq1}). 
 Further, we define $W_n(-\infty)=0$ and $W_n(\infty)=1/\sqrt{n\hat\sigma^2}\sum_{i=1}^ne_i$. Note that $W_n(t)$ is a step function which is constant between consecutive $\hat Y$-order statistics $\hat Y_{(i)}$ and has jumps $e_{[i]}$ there, where $e_{[i]}$ is the residual associated with $i$-th $\hat Y$-order statistic, $\hat Y_{(i)}$. Also, since the model includes the intercept, the following equality holds, $W_n(\infty)=0$, hence rescaling $t$ to $t\in[0,1]$, and using $B_n(t)=W_n(t)-tW_n(1)$, results in identical random process. 



A test statistic, $T$, can then be any function of $W_n(t)$ such that large values give evidence for the lack-of-fit. We considered two commonly used statistics, the KS type statistic
$$T_S=\sup_{t\in \mathds{R}}|W_n(t)|,$$
and the CvM type statistic,
$$T_C=\int_\mathds{R} W_n(t)^2F_n(dt), $$
where $F_n(\cdot)$ is the empirical distribution function of $\hat Y$.

The null distribution of the proposed test statistics is obtained with (random) permutations. 
Define the stochastic process in permutation $k$,  $k=1,...,K$ as
\begin{equation}
\label{permproc}
W_n^k(t)=\frac{1}{\sqrt{n\hat\sigma^2_k}}\sum_{i=1}^ne_i^k I(\hat Y_i^k \leq t), t\in \mathds{R},,
\end{equation}
where
 $$\hat\sigma_k=\sqrt{\frac{1}{n-p}\sum_{i=1}^n (e_i^k)^2} .$$
 (Again, any other consistent estimator of $\sigma$ could be used instead of $\hat\sigma_k$.) With raw data permutation $Y_i^k$ is obtained by randomly permuting $Y_i$. With the permutation of residuals $Y_i^k$ is set to
\begin{equation}
\label{FL}
Y_i^k=\hat Y_i+e_{\pi(i)},
\end{equation}
where $e_{\pi(i)}$ are obtained by randomly permuting $e_i$, $i=1,...,n$. Note that the model is refitted to the permuted outcome and the residuals from the refitted model are used to construct the random process in permutation $k$. 

Since we were mainly considering situations where it was computationally infeasible to calculate the permutation $p$-value based on the whole permutation group, i.e. all $n!$ possible permutations, the $p$-value was calculated  based on random permutations \citep{dwass57,phipsonsmith10,Hemerik2018} using some large value of $K$, where the $p$-values are estimated by using
$$ \frac{1}{K+1}\left(\sum_{k=1}^K I(T^k\geq T)+1\right),$$
 for any test statistic $T$ obtained on $W_n(t)$ and $T_k$ obtained on $W_n^k(t)$, $k=1,\ldots,K$. 
 The usual definition of the permutation $p$-value could be used when performing all possible permutations is feasible (unlikely in the goodness-of-fit context).

To test for the lack-of-fit for a single covariate, say covariate $j$, the tests are modified by ordering the residuals by the respective covariate, i.e using $I(x_{ij}\leq t) $ instead of $I(\hat Y_i \leq t) $ and $I(\hat Y_i^k \leq t) $ in (\ref{eq1}) and (\ref{permproc}), respectively.  It is also straightforward to test for the lack-of-fit for a set of covariates, by using $I(\sum_{j}x_{ij}\hat\beta_j\leq t) $ and $I(\sum_{j}x_{ij}\hat\beta_j^k\leq t) $ in (\ref{eq1}) and (\ref{permproc}), respectively, where the sum is taken only over the defined set of covariates. 
As will be illustrated later, considering the set of specified covariates enables efficient detection of the lack-of-fit due to omission of the interaction effect. Testing the lack-of-fit due to a set of covariates is also important
when the model besides the linear component of the covariate includes also its nonlinear transformation and one wants to test for the adequacy of such assumed nonlinear relation with the outcome. Throughout, we refer to testing the adequacy of the entire model as the full model check, while partial model check refers to testing the adequacy of a model's subset.

\begin{figure}[h!]
\caption{An example of the association between the outcome and the predictor variable (A), permuted outcome using raw data permutation and the predictor variable (B) and permuted outcome using permutation of residuals and the predictor variable (C).}
\label{figex}
\begin{center}
 \resizebox{120mm}{!}{\includegraphics{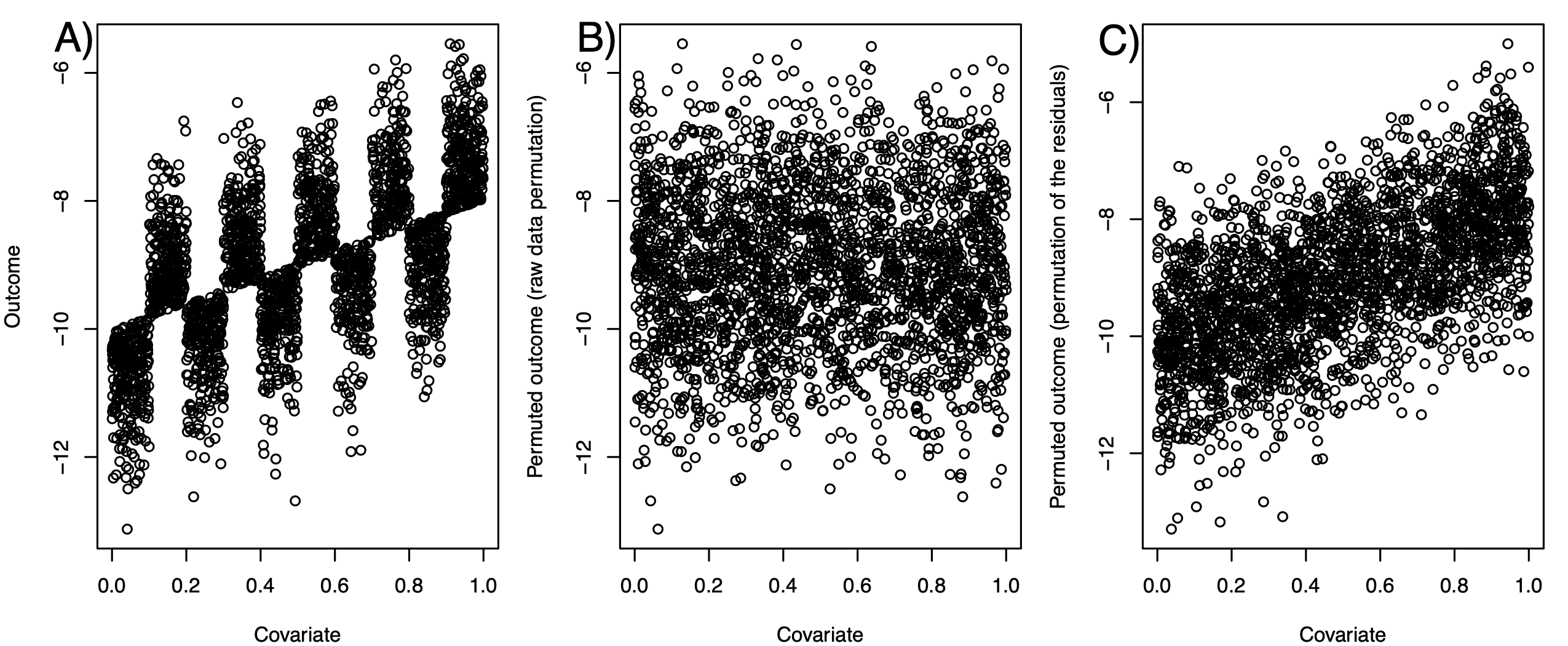}}\end{center}
\end{figure}

The intuitive reason why using permutations is a valid way of approximating the null distribution of the proposed test statistics can be explained as follows. When assuming that the model fits the data well, there should be no obvious pattern when looking at the residuals as a function of the fitted values, while deviations from the model's assumption will appear in such plots as systematic patterns. As an example, assume that the association between the outcome and a single covariate is such as presented in Figure \ref{figex} (A). Clearly, such data could not be generated by a (univariate) linear model. 
However, when using either permutation approach the pattern is no longer present (Figure \ref{figex}, B and C, respectively), obtaining a situation which one observes when the null hypothesis holds.



An important step when constructing the tests is however the standardization of the residuals. Note that by the raw data permutation procedure any positive or negative association between the outcome and the predictor variables is diminished (Figure \ref{figex}, B), which increases the variability of residuals. The random processes constructed under permutations using the unstandardized residuals are therefore more variable and the tests would reject the true null hypothesis less often than they should.  Also, due to increased variance, this permutation approach can suffer from severe loss of power, hence it will not be considered here further (See Supplementary material 2 for the results when using this permutation approach).  

It is necessary to standardize the residuals also when using the permutation of residuals (Figure \ref{figex}, C). In this case the residuals are after permutation and refitting the model less variable (especially obvious with more covariates or when the sample size is small), hence the variability of the random processes constructed after permutations is smaller and the tests would reject the true null hypothesis too often. In fact, the residual variance after the permutations will always be smaller (or equal in the special case where the permuted residuals are equal to the unpermuted residuals) than the residual variance on the data which have not been permuted, which can easily be shown. The effect of standardization is however asymptotically negligible as will be apparent from our theoretical results presented later.


 One could be tempted to estimate the null distribution of the test statistic by constructing the random process directly from the permuted residuals, $e_{\pi(i)}$, without refitting the model. Such test will however be very conservative as this will overestimate the variance of the constructed random process. When the residuals are sorted based on the fitted values, then neighboring residuals are negatively correlated (especially for the residuals close to each other in this ordering). When the residuals are randomly permuted, this negative correlation is preserved, however neighboring residuals are less negatively correlated which overestimates the variance of the random process. When $n$ is large the negative correlation becomes smaller, but the constructed stochastic process depends on a lot of residuals which are negatively correlated, and all these small negative correlations combined are still important. Therefore a new model is fitted after the permutation, to make sure that this negative correlation of residuals that are close with respect to the ordering corresponding to the fitted values is kept.

\subsection{Assumptions}

The theoretical results will be valid under the following set of assumptions.

\begin{itemize}
\item[(A1)] The errors, $\epsilon_i$ are i.i.d. with $E(\epsilon_i)=0$ and $var(\epsilon_i)=\sigma^2<\infty$, $i=1,\ldots,n$. We further assume that $\lim_{n\rightarrow \infty}\frac{1}{n}\sum_{i=1}^nE|c_{ij}\epsilon_i|=0$ for any constants $0<c_{ij}<\infty$, $i=1,\ldots,n$, $j=1,\ldots,k$.
\item[(A2)] $0<\lim_{n\rightarrow\infty}||\frac{1}{n}\mathbf{x}^T\mathbf{x}||_2<\infty$.
\item[(A3)] $ E\left[   \epsilon_i I( \mathbf{x}_i^T\bm{\gamma}_1 \leq x_1   ) -  \epsilon_i I( \mathbf{x}_i^T\bm{\gamma}_2 \leq x_2   )  \right]^2\rightarrow 0$, as $\bm\gamma_1\rightarrow \bm\gamma_2$ and $x_1\rightarrow x_2$, for all $\bm\gamma_2$ sufficiently close to $\bm\theta$.
\item[(A4)] When $n\rightarrow\infty$, the functions
$$h(t)=(h_1(t),\ldots,h_{p}(t)) $$
and
$$k(t)=(k_1(t),\ldots,k_{p}(t)), $$
with $h_j(t)=\sum_{i=1}^n x_{ij}I( \mathbf{x}_i^T\bm{\theta} \leq t   ) $ and $k_j(t)=\sum_{i=1}^n x_{\pi(i)j}I( \mathbf{x}_i^T\bm{\theta} \leq t   )$, $j=1,\ldots,p$, respectively, are continuous for each $t\in \mathds{R}$.
\end{itemize}

Assumption (A1) is standard and assures that the weak law of large numbers and the central limit theorem can be applied to sequences $n^{-1}\sum_{i=1}^n c_i \epsilon_i$ and $n^{-1/2}\sum_{i=1}^n c_i \epsilon_i$, respectively, for any $0<c_{i}<\infty$. Note that the strong law of large numbers can be applied to a sequence $n^{-1}\sum_{i=1}^n  \epsilon_i$.

Assumption (A2) implies that $\frac{1}{n}\mathbf{x}^T\mathbf{x}$ and $(\frac{1}{n}\mathbf{x}^T\mathbf{x})^{-1}$ are finite positive-definite matrices when $n\rightarrow\infty$. It also implies that $0<\lim_{n\rightarrow\infty} ||(\frac{1}{n}\mathbf{x}^T\mathbf{x})^{-1}||_2<\infty $.

Assumption (A3) is equivalent to condition (v) in the functional central limit theorem in \citet[p. 53]{Pollard90} and assures stochastic equicontinuity of $W_n(t)$ and $W_n^k(t)$. Assumption  (A4) together with assumption (A3) guarantees that in the limit $W_n(t)$ and $W_n^k(t)$ will have continuous sample paths.

The results presented herein could be adapted also for the random design matrix, but at a cost of some additional distributional assumptions and a large computational inconvenience. Also, while we consider here only the OLS estimator, the results would hold for any consistent estimator of $\bm\theta$.

\subsection{Convergence of the constructed stochastic processes under the null hypothesis}

Here we state the main theoretical result of this paper. That is, that under the hypothesis that the data were generated by model (\ref{eq1.1}) permutation of the residuals  yields a valid approximation of the constructed random process.

\begin{theorem}
With permutation of residuals, under $H_0$ and (A1), (A2), (A3) and (A4), we have with probability 1, that the process $\frac{\hat\sigma}{\sigma}W_n(t)$ and, conditionally on the data, the process $\frac{\hat\sigma_k}{\sigma}W_n^k(t)$  converge in distribution to the same zero mean Gaussian process $G_\infty$ in the Skorokhod-space $D[-\infty,\infty]$,  where the covariance function of $G_\infty$ is


$$K(t,s)=\lim_{n\rightarrow\infty}\frac{1}{n}\left( \sum_{i=1}^n  I( \mathbf{x}_i^T\theta \leq t \wedge s  ) - h(t)(\mathbf{x}^T\mathbf{x})^{-1} h(s)^T \right).$$
where
$$h(z)= \sum_{i=1}^n \mathbf{x}_{i}^TI(\mathbf{x}_{i}^T\theta \leq z).$$

\end{theorem}
\begin{remark}
Let $g$ be some continuous function from the Skorokhod space $D[-\infty,\infty]$ to the real line. From the consistency of $\hat\sigma$ and $\hat\sigma_k$ (see Supplementary material 1), the continuous mapping theorem and Slutsky's theorem it then follows that $\frac{\sigma}{\hat\sigma}g( \frac{\hat\sigma}{\sigma}W_n(t) )$ and, conditionally on the data, $\frac{\sigma}{\hat\sigma_k}g( \frac{\hat\sigma_k}{\sigma}W_n^k(t) )$ have the same limiting distribution.
\end{remark}
\begin{remark}
The result also applies when one performs partial model checks,  where the covariance function of the limiting random process is then,
$$K(t,s)=\lim_{n\rightarrow\infty}\frac{1}{n}\left( \sum_{i=1}^n  I( \sum_{j}x_{ij}\beta_j \leq t \wedge s  ) - h(t)(\mathbf{x}^T\mathbf{x})^{-1} h(s)^T \right),$$
where
$$h(z)= \sum_{i=1}^n \mathbf{x}_{i}^TI(\sum_{j}x_{ij}\beta_j \leq z),$$
and the sum $\sum_{j}$ is taken only over the defined set of covariates.
\end{remark}
\begin{remark}
When specialized to our setting, the tests proposed by \citet{Su91,Lin02,Stute98a} can be seen as special cases of the proposed procedure (see Supplementary material 1) hence the result applies  also for their tests. The result applies also when the residuals are not standardized (see Supplementary material 1) where the covariance function is then
$$K(t,s)=\lim_{n\rightarrow\infty}\frac{\sigma^2}{n}\left( \sum_{i=1}^n  I( \mathbf{x}_i^T\theta \leq t \wedge s  ) - h(t)(\mathbf{x}^T\mathbf{x})^{-1} h(s)^T \right).$$
The result also applies when the refitted residuals are ordered by the original fitted values which implies that in the approach proposed by \citet{Su91,Lin02} $e_iZ_i$ in (\ref{sw}) can be replaced by the $i$th permuted residual (note that in this case it is not necessary to refit the model for each $m=1,\ldots,M$). The result also applies to the tests based on partial model checks \citep{Christensen15,Hattab18}, implying that the proposed procedure could also be used as an alternative of the Monte-Carlo procedure proposed by \citet{Hattab18}.
\end{remark}

\subsection{Convergence of the constructed stochastic process under the alternative hypothesis}

The alternative hypothesis which we are interested in is that there does not exist any $p$-vector $\mathbf{b}$ such that $\nu(\mathbf{x})=E(Y|\mathbf{x})=\mathbf{x}^T\mathbf{b}$. Note that under $H_1$, $\bm{\hat\theta}$ converges to a constant vector $\bm{\theta}_*$ in probability as $n\rightarrow \infty$.

\begin{proposition}
\label{lemaH1}
Let $g$ be some continuous function from the Skorokhod space $D[-\infty,\infty]$ to the real line. Under $H_1$, $\frac{1}{\sqrt{n}}g(W_n(t))$ converges in probability towards some non-zero constant $c\neq 0$.
\end{proposition}
\begin{remark}
The result can be applied directly to establish that $\frac{1}{\sqrt{n}}T_C$ and also $\frac{1}{\sqrt{n}}T_{C^*}= \frac{1}{\sqrt{n}}\int_\mathds{R} W_n(t)^2dt$ converge in probability towards some positive constants, $c_1> 0$ and $c_2>0$, respectively.
\end{remark}



\begin{proposition}
\label{col1}
Under $H_1$, 
$\frac{1}{\sqrt{n}}T_S$ converges in probability towards some positive constant $c_1>0$.
\end{proposition}

\section{Simulation results for the full model check}

To check the size of the tests and to illustrate their power for several alternatives we conducted a simulation study with two predictor variables and an intercept term. The fitted model is always,

$$y=\beta_0+\beta_1x_1+\beta_2x_2+\epsilon. $$

The predictor variables were simulated independently from the uniform $[0,1]$ distribution for $n=50$, $100$, $200$, $500$ and $1000$ subjects, if not stated otherwise. In each step of the simulation $K=1000$ permutations were performed. Each step of the simulation was repeated $10,000$ times (simulation margins of error are $\pm 0.002$, $\pm 0.004$ and $\pm 0.006$ for $\alpha=0.01$, 0.05 and $0.1$, respectively).

All the analyses were performed with R language for statistical computing, version 3.0.3 \citep{R}. R package \verb"gof" was used for the tests proposed by \citet{Lin02} while our proposed tests are available through R package \verb"gofLM" (available upon request from the authors).

\subsection*{Size of the test}

The size of the tests was verified with a simulated example where the outcome variable was simulated from,

$$y_i=-0.1+\beta_1 x_{i1}+ \beta_2 x_{i2}+ \epsilon_i.$$

\subsubsection*{Normal homoscedastic random errors}
Here $ \epsilon_i\sim N(0,\sigma^2)$, where $N(\cdot)$ denotes a normal distribution. Different values of $\beta_1=0$, $0.25$, $0.50$, $0.75$ and $1$, $\beta_2=0$ and $0.25$ and $\sigma^2=0.01$, $0.10$, $0.25$, $0.50$, $0.75$, $1$ and $2$ were considered. Sample size ranged from $n=5$, $7$, $10$, $50$, $100$, $200$, $500$ and $1000$ subjects. While performing a goodness-of-fit test with $n<50$ is in practice highly questionable, $n<50$ is included here for the sake of illustration.



\begin{figure}[h!]
\caption{Distribution of the $p$-values, based on 10,000 simulated data sets under the null hypothesis using normal random errors, $\beta_1=\beta_2=0.25$, $\sigma^2=0.25$. Columns: different tests, rows: $n$.}
\label{p.nul}
\center{\resizebox{120mm}{!}{\includegraphics{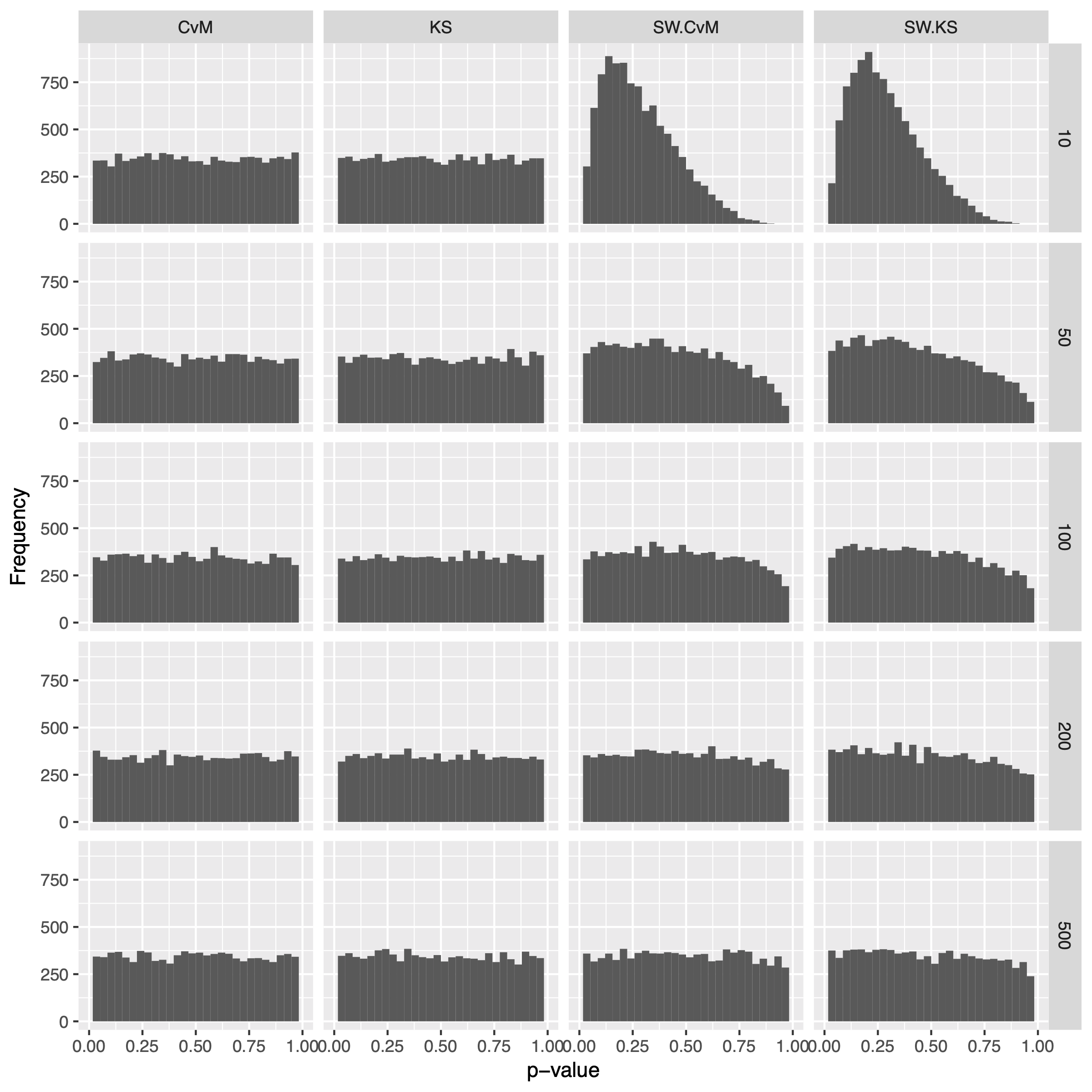}}}
\end{figure}

The distribution of $p$-values for the situation with $\beta_1=\beta_2=0.25$ and $\sigma^2=0.25$ and some selected sample sizes using the tests proposed by \citet{Lin02} (SW's test) and our proposed tests using KS and CvM type test statistics are shown in Figure \ref{p.nul}. It is obvious that our proposed tests perform well, while the SW's tests do not perform well in terms of size, since the distribution of the $p$-values is not entirely uniform event with $n=500$. It is shown in the Supplementary material that with a very small sample size our proposed tests were slightly conservative, however with as few as 10 subjects, the tests attained the nominal level. Hence we can conclude that the convergence rate of the proposed tests is very fast, which agrees also with Figure \ref{p.nul}. 
In agreement with Figure \ref{p.nul}, the rejection rates of the SW's tests are outside the simulation margin of error for small to moderate sample sizes (Supplementary material 2). 



\subsubsection*{Non-normal homoscedastic random errors}

The size of the proposed tests with non-normal random errors was checked with the example where $\epsilon_i=\epsilon_i^1-as$, $ \epsilon_i^1\sim \gamma(a,s)$ and $a$ and $s$ are the shape and the scale parameter of the gamma distribution ($\gamma(\cdot)$), respectively using $\beta_1=\beta_2=0.25$. The scale parameter was set to $s=1$, while different values of the shape parameter were considered, $a=1$, $2$, $5$, $7$, $10$ and $15$. Note that larger values of $a$ mean that distribution of the error term was more symmetric but also more variable. The sample sizes considered here were, $n=5$, $7$, $10$, $15$, $20$, $30$, $40$, $50$, $100$, $200$, $500$ and $1000$.




With non-normal random errors the proposed tests were also slightly conservative with a very small sample size ($n=5$) but with a larger sample size the distribution of $p$-values was uniform and the tests attained the nominal level (Supplementary material 2). The distribution of $p$-values obtained by SW's tests was not uniform even with $n=500$ and the tests did not attain the nominal level (Supplementary material 2).  

\subsubsection*{Normal heteroscedastic random errors}

Here we illustrate the performance of the proposed tests in the presence of heteroscedasticity, with $\epsilon_i\sim N(0,\sigma^2_i)$, where $\sigma^2_i=(1+\theta x_{i1})^2$. Different values of $\theta=-0.5$, $-0.2$, $-0.1$, $-0.05$, $0$, $0.05$, $0.10$, $0.20$, $0.30$ and $0.50$ were considered. The other parameters were: $\beta_1=0$, $0.25$, $0.50$, $1$ and $2.5$, $\beta_2=0.25$.





As expected, the tests were liberal in the presence of heteroscedasticity.  The rejection rates were however only slightly larger than the nominal levels even when there was strong heteroscedasticity (large $|\theta|$), especially when using the CvM type test statistic (Supplementary material 2).


\subsection*{Omission of a quadratic effect}

Here we illustrate the power of the proposed test to detect the lack-of-fit due to omitting the quadratic effect. The outcome variable was simulated from,

$$y_i=-0.1+\beta_1x_{i1}+0.25x_{i2}+\beta_3 x_{i1}^2  +\epsilon_i,$$
$ \epsilon_i\sim N(0,0.1)$. Different values of $\beta_1=0$, $0.25$, $0.50$ and $1$ and $\beta_3=0$, $0.25$, $0.50$, $0.75$, $1$ were considered.

Results for different nominal test size $\alpha$, for different effect and sample sizes for the test using CvM test statistic and permutation of residuals are reported in Figure \ref{fig1both}, lower panels. The results for the other tests and settings are reported in Supplementary material 2.

\begin{figure}[h!]
\caption{Fraction rejected for the goodness-of-fit test; based on 10,000 simulated data sets under the null hypothesis ($\beta_3=0$) and under different alternatives (omission of a quadratic term, omission of an interaction term; $\beta_3>0$), $\beta_1=0.25$. Cramer-von Mieses type test statistic. Columns: $\alpha$.}
\label{fig1both}
\center{\resizebox{120mm}{!}{\includegraphics{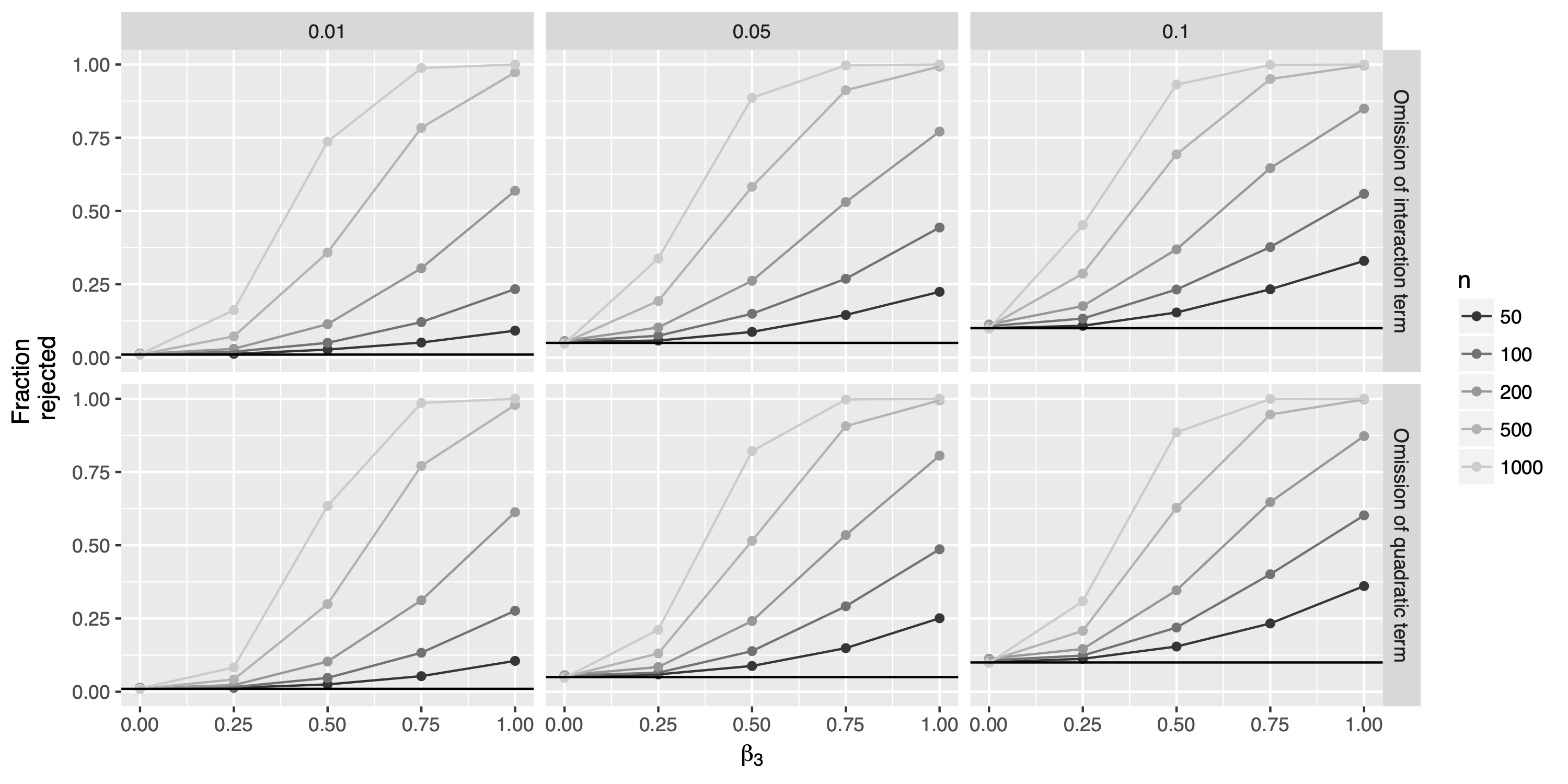}}}
\end{figure}

Judging from Figure \ref{fig1both} (upper panels) the test  performs well. As expected, the test was more powerful when $\beta_1$, the effect size ($\beta_3$) and sample size was larger. The test using CvM type test statistic was more powerful than the test based on KS type test statistic (Supplementary material 2).

\subsection*{Omission of an interaction term}

In this example we illustrate the power of the proposed test to detect the lack-of-fit due to ignoring the interaction effect. The outcome variable was simulated from,

$$y_i=-0.1+\beta_1x_{i1}+0.25x_{i2}+\beta_3 x_{i1}x_{i2}  +\epsilon_i,$$
$ \epsilon_i\sim N(0,0.1)$. Different values of $\beta_1=0$, $0.25$, $0.50$, $1$ and $\beta_3=0$, $0.25$, $0.50$, $0.75$ and $1$ were considered.


Results for different nominal test size $\alpha$ for test using CvM test statistic and permutation of residuals are reported in Figure \ref{fig1both}, upper panels. The results for the other tests and settings are reported in Supplementary material 2.

The test was powerful to detect the lack-of-fit due to omission of the interaction term. In this example the test based on CvM type test statistic was again more powerful. 


\section{Simulations results for partial model checks}

Simulations were performed also to check the properties of the test where one tests for the lack-of-fit due to a specific covariate or a set of covariates. Here there were 5 variables and the intercept term, and the fitted model is,

$$y=\beta_0+\beta_1x_1+\beta_2x_2+\beta_3x_3+\beta_4x_4+\beta_5x_5+\epsilon. $$

 The predictor variables were simulated independently from the uniform $[0,1]$ distribution for $n=50$, $100$, $200$, $500$ and $1000$ subjects. In each step of the simulation $K=1000$ permutations were performed. Each step of the simulation was repeated 10,000 times. Here we only report the results for $n=100$ and for the tests based on the CvM type test statistic, the other results are reported in Supplementary material 2.

\subsection*{Example I}

In this example the size and power of the tests against the alternative where one of the covariates, $x_1$, has a nonlinear association is shown. The outcome variable was simulated from,

$$y_i=-0.1+0.25x_{i1}+0.25x_{i2}+0.25x_{i3}+0.25x_{i4}+0.25x_{i5} +\beta_6 x_{i1}^2 +\epsilon_i,$$
$ \epsilon_i\sim N(0,0.1)$. Different values of $\beta_6=0$, $0.25$, $0.50$, $0.75$, $1$ were considered.

\begin{figure}[h!]
\caption{Fraction rejected for the goodness-of-fit tests; based on 10,000 simulated data sets with the quadratic effect (example I) and the interaction effect (example II). Example I: full model check (global) and different tests targeting at the correct variable, $x_1$ (correct) and the wrong variable, $x_2$ (wrong); example II: full model check (global) and different tests targeting at the set of covariates (targeted at the set of correct variables, $x_1$ and $x_2$ (2 correct), two wrong variables, $x_3$ and $x_4$ (2 wrong), one correct and one wrong variable, $x_1$ and $x_3$ (1 correct 1 wrong) and 3 wrong variables, $x_3$, $x_4$ and $x_5$ (3 wrong)). The tests are based on the Cramer-von Mieses type test statistic (CvM). Columns: $\alpha$.}
\label{figpostfin}
\center{\resizebox{120mm}{!}{\includegraphics{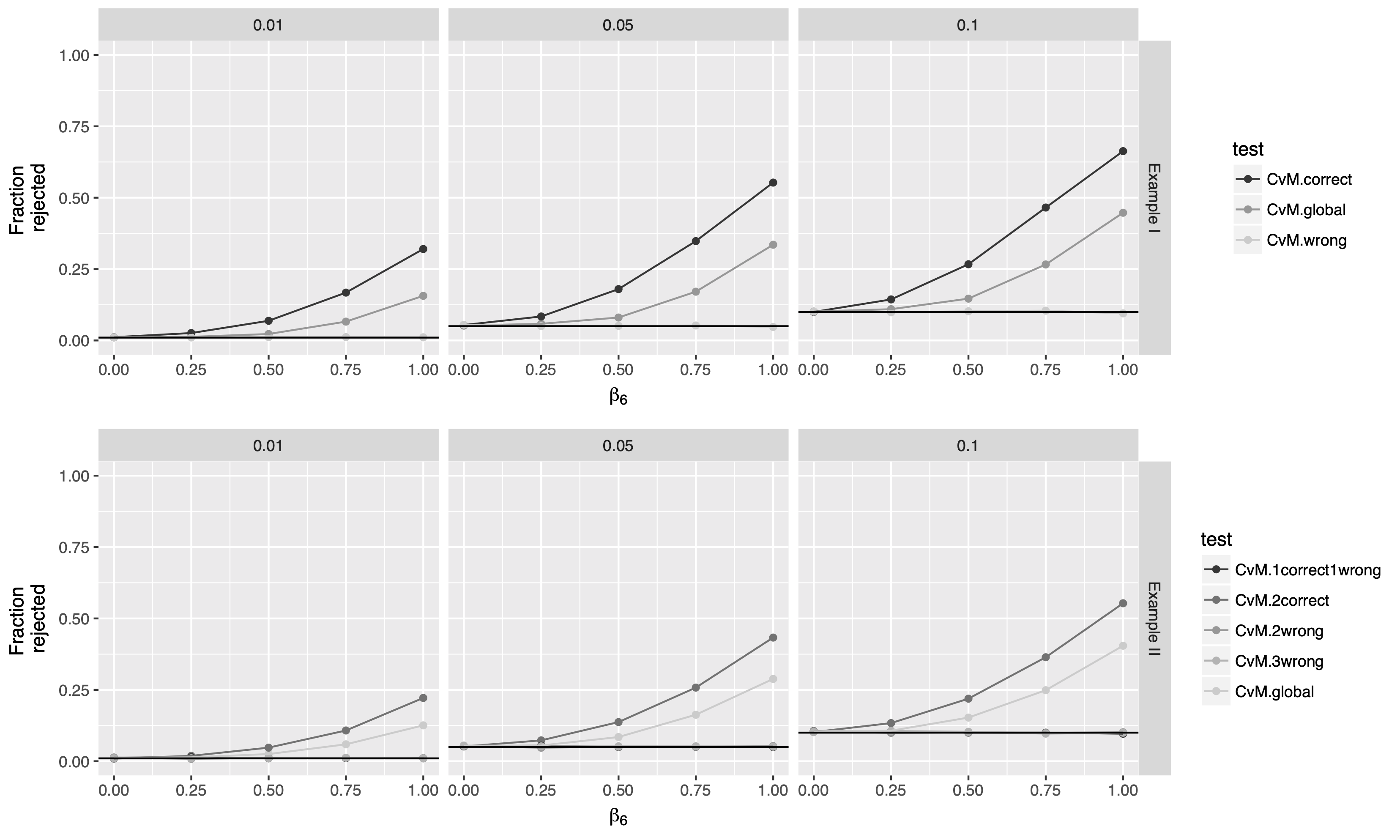}}}
\end{figure}

The fraction of rejected null hypotheses for the full model check as well as for the partial model check targeting only at the variable which has a nonlinear association ($x_1$), and the partial model check targeting only the variable which has a linear association with the outcome ($x_2$), are reported in Figure \ref{figpostfin} (upper panels).



The test targeting at $x_1$ had good power and the rejection rates were, as expected, larger than the rejection rates of the full model check. The test targeting at $x_2$ attained the nominal level. The tests based on the CvM type test statistic were more powerful (Supplementary material 2).


\subsection*{Example II}

In this example the size and power against the alternative where there is an interaction between covariates $x_1$ and $x_2$ is shown. The outcome variable was simulated from,

$$y_i=-0.1+0.25x_{i1}+0.25x_{i2}+0.25x_{i3}+0.25x_{i4}+0.25x_{i5}+\beta_6 x_{i1}x_{i2}  +\epsilon_i,$$
$ \epsilon_i\sim N(0,0.1)$. Different values of $\beta_6=0$, $0.25$, $0.50$, $0.75$ and $1$ were considered.

The fraction of rejected null hypotheses for the full model check as well as for the tests targeting at $x_1$ and $x_2$, and the wrong variable, i.e. the variable $x_3$, are reported in Supplementary material 2.

While the tests for the full model check were powerful when $\beta_6> 0$, all tests targeting at a specific covariate had no power in this example as the rejection rates were the same as the nominal level. As expected, targeting at only one variable has in this case no power as the lack-of-fit is due to a set of covariates. Therefore, when one wants to detect the lack-of-fit due to omission of the interaction term, a set of suspect variables needs to be defined. The results when the test targeted at the set of correct variables, i.e. variables $x_1$ and $x_2$, one correct and one wrong variable ($x_1$ and $x_3$), two wrong variables ($x_3$ and $x_4$) and three wrong variables ($x_3$, $x_4$ and $x_5$)  are reported in Figure  \ref{figpostfin} (lower panels). 

The test targeting at the correct set of covariates had good power, while the tests targeting at the subsets of variables which contain also the variables which only have the main effect with the outcome attain the nominal level. Tests based on CvM type test statistic were more powerful (Supplementary material 2).







\section{Applications}
Here we apply the proposed tests to two publicly available data sets.

\subsection*{Application 1}
\begin{figure}[h!]
\caption{The constructed random process (red curves) and 1000 randomly selected random processes which are obtained after using permutation of residuals (gray curves) for the fetal mandible length data. The $p$-values are for the CvM type test statistic. Panel (A) refers to the model which includes only the effect of age, while panel (B) is the model which includes also the square of age.}
\label{figReal}\center{
\center{\resizebox{120mm}{!}{\includegraphics{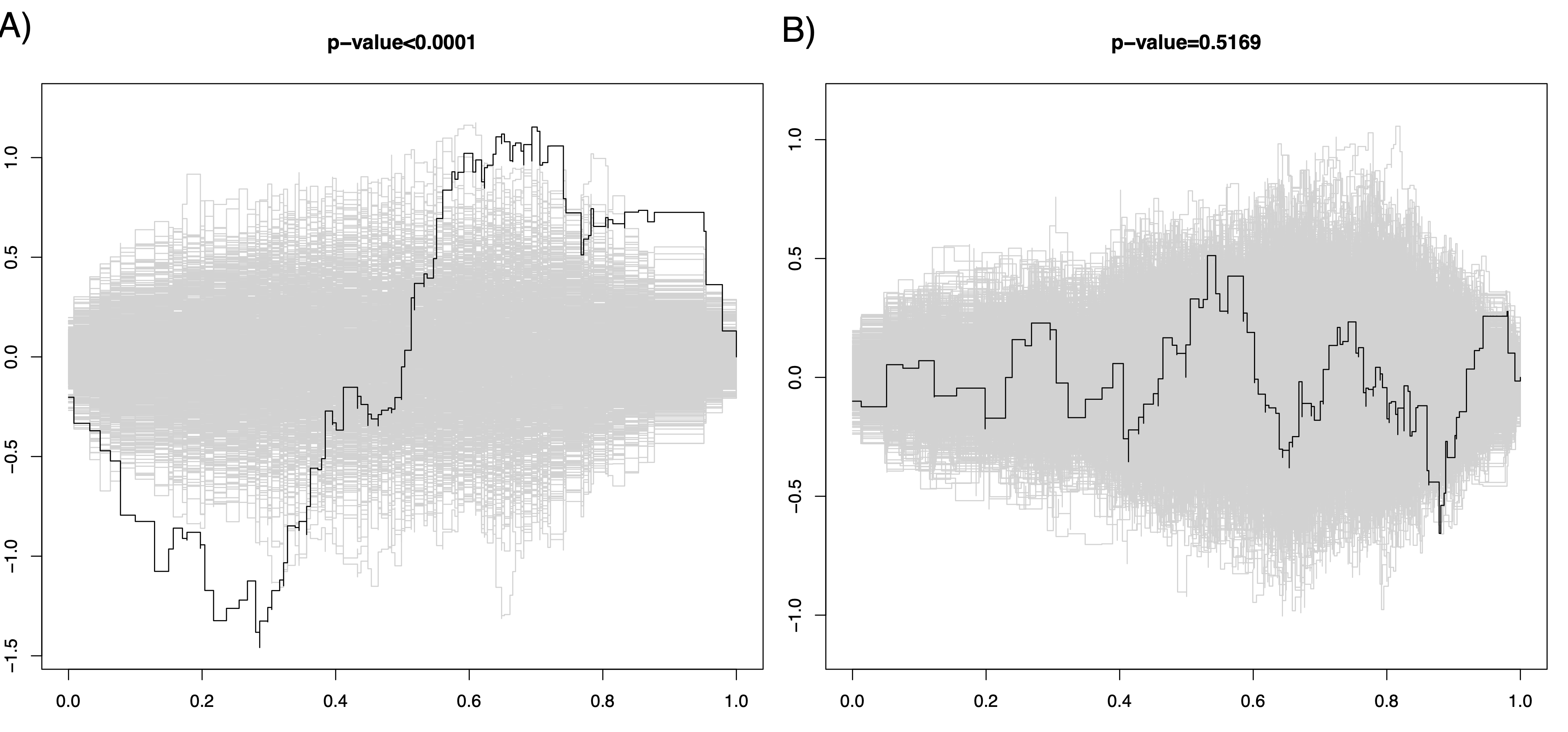}}}}
\end{figure}

The proposed tests were applied to data from a study of fetal mandible length by \citet{Chitty93}. The data comprises measurements of mandible length and gestational age in 158 fetuses. The data are publicly available through R package \verb"lmtest". The data were log transformed and the proposed goodness-of-fit tests, using 10,000 permutations, all returned a $p$-value of $<0.0001$ and the null hypothesis was rejected. When looking at the plot of the constructed random process (Figure \ref{figReal}, panel (A), red curve) we see that there is a long succession of negative residuals, followed by a long succession of the positive residuals, which suggests that there is a quadratic association. When the goodness-of-fit tests were applied to the model which included also the quadratic term for gestational age the $p$-values were 0.426 and 0.517 for the test using KS and CvM type test statistic, respectively and the null hypothesis was not rejected at $\alpha=0.05$. In this case the constructed random process shows no systematic pattern (Figure \ref{figReal}, panel (B)) and is similar as the random processes which are obtained after using permutation of residuals (gray lines in Figure \ref{figReal} are 1000 randomly selected processes which are obtained after using permutation of residuals).



\subsection*{Application 2}
\begin{figure}[h!]
\caption{The constructed random process (red curves) and 1000 randomly selected random processes which are obtained after using permutation of residuals (gray curves) for the Draper and Smith data. The $p$-values are for the CvM type test statistic. Panel (A) refers to the model which includes only the main effects. Panels (B) and (C) are tests targeting at $x_1$ and $x_2$ for the model which includes only the main effects. Panel (D) refers to the model which includes also $x_1^2$. Panel (E) refers to the test targeting at $x_1$ and $x_1^2$ and panel (F) is for a test targeting at $x_2$.}
\label{figReal2}\center{
\center{\resizebox{120mm}{!}{\includegraphics{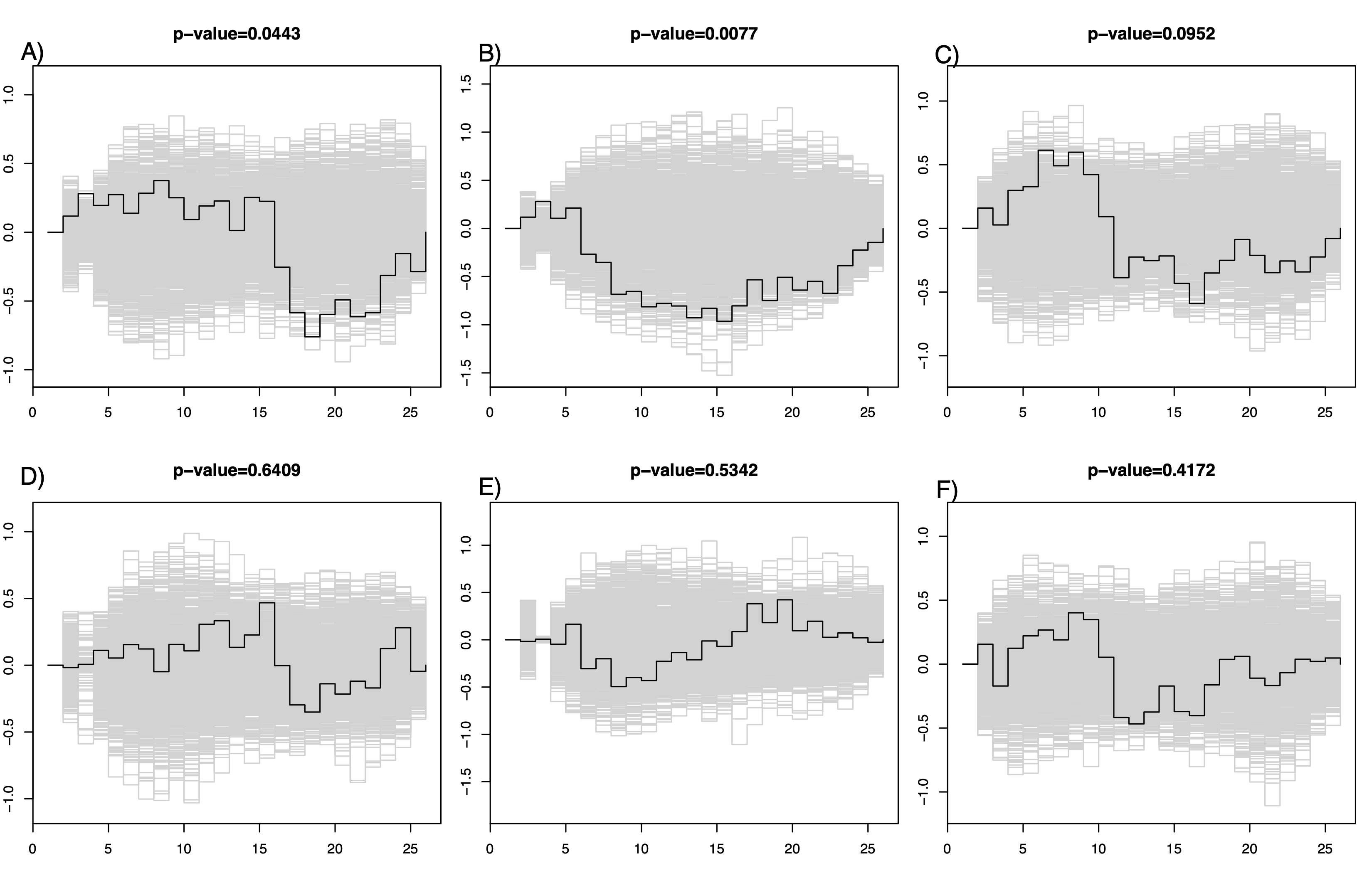}}}}
\end{figure}

The tests were applied also to a dataset which was taken from \citet{draper81}. This dataset consists of two covariates (the operating days per month, $x_1$, and the average atmospheric temperature, $x_2$) and the pounds of steam used monthly as the outcome for 25 subjects. \citet{draper81} used a simple additive model and concluded, based on the plot of residuals versus predicted values that the data fitted the model well. Our tests which use permutation of residuals returned $p$-values 0.042 and 0.044 for KS and CvM type test statistic, respectively and the null hypothesis was rejected at $\alpha=0.05$ (Figure \ref{figReal2}, panel (A)). The tests using CvM type test statistic targeting at $x_1$ and $x_2$ returned $p$-values equal to 0.008 and 0.096, respectively. From the plots of the constructed random process we can conclude that the lack-of-fit due to covariate $x_1$ is due to omission of the quadratic term (Figure \ref{figReal2}, panel (B), red curve). When the quadratic term for $x_1$ was included in the model, the $p$-values for the full model check were 0.567 and 0.641 for KS and CvM type test statistic, respectively and the null hypothesis was not rejected at $\alpha=0.05$ (Figure \ref{figReal2}, panel (D)). The test using CvM type test statistic targeting at $x_1$ and its square, $x_1^2$, returned a $p$-value of 0.534 and the test targeting at $x_2$ returned a $p$-value of 0.417. The plots of the constructed process when targeting at $x_1$ and $x_1^2$ reveal that the goodness-of-fit improved dramatically by including $x_1^2$ in the model (Figure \ref{figReal2}, panel (E), red curve). This had an effect also on the random process which targets at $x_2$, where after inclusion of $x_1^2$ the random process is more in line with what one observes under the null hypothesis (Figure \ref{figReal2}, panel (F), red curve).

 \section{Discussion}

In practice, graphical procedures are mainly used to evaluate goodness-of-fit of the linear regression model. While these procedures are useful, they can be to a large extent subjective. They can be formalized by forming a stochastic process based on residuals. Several goodness-of-fit tests based on such random processes have been proposed. They are however very computationally expensive \citep{Stute98a}, target  a certain fixed alternative \citep{Lin02,Aldo12}, rely on the asymptotic independence of the residuals \citep{Fan01,Christensen15} or rely on some asymptotic properties of the constructed random process \citep{Su91,Lin02} and therefore can yield poor results with small sample size especially with non-normal random errors. We propose goodness-of-fit tests which do not depend on the distribution of the errors and have correct size also with a small sample size. This is achieved by standardizing the residuals and using permutations to obtain the $p$-values. We prove that when using the permutation of residuals, under homoscedastic random errors, the tests are asymptotically consistent under the null and the alternative hypotheses. 

We also show how the proposed tests can be used to detect the lack-of-fit which occurs only due to a specific covariate or a set of covariates. This is achieved by using different ordering of the residuals. To detect the lack-of-fit only due to one covariate then the residuals are ordered by this covariate, while when trying to detect the lack-of-fit of the model's subset, the residuals are ordered by the predictive values obtained from the respective model's subset. Other orderings that were proposed, e.g. the moving sums and moving averages proposed by \citet{Lin02}, could also easily be implemented, but were not pursued in this study.  Also, the standardization and ordering techniques used by \citet{Hattab18} could all easily be adopted also for the cumulative sum processes that are considered here. The proposed procedure could also be used in tests based on partial sums of residuals as an alternative to the Monte-Carlo procedure proposed by \citet{Hattab18}.

A possible issue is the performance of the proposed tests in the presence of heteroscedastic random errors. We show with extensive simulations that even in the presence of strong heteroscedasticity the type I error is only slightly inflated, especially when using the CvM type test statistic. However, the wild bootstrap proposed by \citet{Stute98a} or the approximation proposed by \citet{Lin02} are in such settings superior to our tests, as they have correct size, at least asymptotically. 

Similarly as the tests proposed by \citet{Stute98a} and \citet{Lin02}, our tests could be extended to the other generalized linear models (GLM), but the assumption of homoscedastic errors could be problematic for some GLMs. 
We are currently exploring the possibility of extending this idea to linear mixed effects models.

One potential disadvantage of using permutation testing could be a large computational burden. In our case this was not an issue, as it is possible to implement the tests in a very efficient manner. 
E.g. the calculation of $p$-values for our real data example from a study of fetal mandible length using 10,000 permutations took 
1.2 seconds 
on an ordinary personal computer. A typical simulated data example with two covariates, 1,000 samples and 10,000 permutations took 
2.6 seconds. 
In comparison, the tests proposed by \citet{Lin02} using the \verb"gof" package with 10,000 simulations took 28.6 seconds for the same example.   





\newpage
\begin{center}
{\large\bf SUPPLEMENTARY MATERIAL}
\end{center}

\begin{description}


\item[Title: Supplementary material 1.] Theoretical results. In this supplementary material theoretical results are proved. (pdf)
\item[Title: Supplementary material 2.] Additional simulation results. In this supplementary material additional simulation results are shown. These include, among others, also the simulation settings considered in \citet{Christensen15,Stute98a}, where the goal was to compare the performance of our proposed tests with the existing goodness-of-fit tests under the simulation design considered by the respective authors. (pdf)

\begin{itemize}
\item \href{https://www.dropbox.com/s/0j87aqvryix6v1f/supplementaryMaterial1.png?dl=0}{Supplementary material 1.(link)}
\item \href{https://www.dropbox.com/s/kuqejm037hfxsch/supplementaryMaterial2.png?dl=0}{Supplementary material 2.(link)}
\end{itemize}

\end{description}

\section*{Acknowledgements}
This work was supported by the Slovenian Research Agency (Predicting rare events more accurately, N1-0035; Methodology for data analysis in medical sciences, P3-0154).

\bibliographystyle{dinat}
\bibliography{gof-LM}

\end{document}